\documentclass[11pt,psfig]{article}
\usepackage{amssymb,amsmath,amscd}
\newtheorem{lem}{Lemma}
\newtheorem{conj}{Conjecture}
\newtheorem{thm}{Theorem}
\newtheorem{cor}{Corollary}

\newcommand{\commentout}[1]{}

\newcommand{\arr}[1]{\left( \begin{array}{cllcr} #1 \end{array} \right)}

\begin{document}
\title{Eigenvectors and Reconstruction}
\author{Hongyu He \\
Department of Mathematics \& Statistics \\
Louisiana State University \\
email: hongyu@math.lsu.edu \\
} 
\date{}
\maketitle 
\abstract{In this paper, we study the simple
eigenvectors of two hypomorphic matrices using linear algebra. We also give a new proof of a result of Godsil-McKay.  }
\section{Introduction}
We start by fixing some notations.
Let $A$ be a $n \times n$ real symmetric matrix. Let $A_i$ be the
 matrice obtaining by deleting the $i$-th row and $i$-th column of $A$. We say that two symmetric matrices $A$ and $B$ are hypomorphic if $B_i$ can be obtained by permuting the rows and columns of $A_i$ simultaneously. Let $\Sigma$ be the set of permutations. We write $B=\Sigma(A)$.\\
\\
 If $M$ is a symmetric real matrix, then the eigenvalues of $M$ are real. We write
$$eigen(M)=(\lambda_1(M) \geq \lambda_2(M) \geq \ldots \geq \lambda_n(M)).$$
If $\alpha$ is an eigenvalue of $M$, we denote the corresponding eigenspace by $eigen_{\alpha}(M)$. 
Let $\mathbf 1$ be the $n$-dimensional vector $(1,1,\ldots,
1)$.  Put $\mathbf J =\mathbf
1^t \mathbf 1$. 

\begin{thm}[~\cite{he}]~\label{main}
Let $B$ and $A$ be two real $ n \times n$ symmetric matrices.
Let $\Sigma$ be a hypomorphism such that $B=\Sigma(A)$. Let $t$ be a real number. Then there exists an open interval $T$ such that for $t \in T$ we have
\begin{enumerate}
\item $\lambda_n(A+ t \mathbf J)=\lambda_n(B+t \mathbf J)$;
\item
$eigen_{\lambda_n}(A+ t \mathbf J)$ and $eigen_{\lambda_n}(B+t \mathbf J)$ are both one dimensional;
\item
$eigen_{\lambda_n}(A+ t \mathbf J)=eigen_{\lambda_n}(B+t \mathbf J).$
\end{enumerate}
\end{thm}
As proved in ~\cite{he}, our result implies Tutte's theorem which says that $eigen(A+ tJ)=eigen(B+t J)$.\\
\\
In this paper, we shall study the eigenvectors of $A$ and $B$. We first prove that the squares of the entries of simple unit eigenvectors of $A$ can be reconstructed as functions of $eiegn(A)$ and $eigen(A_i)$. This yields a proof of a Theorem of Godsil-McKay. We also study how the eigenvectors of $A$ change after a purturbation of a rank $1$ symmetric matrices. Combined with Theorem ~\ref{main}, we prove another result of  Godsil-McKay which states that the simple eigenvectors that are perpenticular to $\bold 1$ are reconstructible. We further show that the orthogonal projection of $\bold 1$ onto higher dimensional eigenspaces is reconstructible.\\
\\
Our investigation indicates that the following conjecture could be true.
\begin{conj} Let $A$ be a real $ n \times n$ symmetric matrix. Then there exists a subgroup $G(A) \subseteq O(n)$ such that a real symmetric matrix $B$ satisfies the properties that $eigen(B)=eigen(A)$ and $eigen(B_i)=eigen(A_i)$ for each $i$ if and only if $B = U A U^t$ for some $U \in G(A)$. 
\end{conj}
This conjecture is clearly true if $rank(A)=1$. For $rank(A)=1$, the group $G(A)$ can be chosen as
$\mathbb Z_2^n$, all in the form of diagonal matrices. In some other cases, $G(A)$ can be a subgroup of the permutation group $S_n$.
\begin{conj}
The group $G(A)$ can be chosen to be a twisted product of a subgroup of $S_n$ with $\mathbb Z_2^n$.
\end{conj}
Clearly, this conjecture implies the reconstruction conjecture.

\section{Reconstruction of Square Functions}

\begin{thm}Let $A$  be a $n \times n$ real symmetric matrix. Let $(\lambda_1 \geq \lambda_2 \geq \cdots \geq \lambda_n)$ be the eigenvalues of $A$.
Suppose $\lambda_i$ is a simple eigenvalue of $A$. Let $\mathbf p_i=(p_{1,i}, p_{2,i}, \ldots, p_{n,i})^t$ be a unit vector in $eigen_{\lambda_i}(A)$. Then for every $m$,
$p_{m, i}^2$ can be expressed as a function of $eigen(A)$ and $eigen(A_m)$.
\end{thm}
Proof: Let $\lambda_i$ be a simple eigenvalue of $A$. Let $\mathbf p_i=(p_{1,i}, p_{2,i}, \ldots, p_{n,i})^t$ be a unit vector in $eigen_{\lambda_i}(A)$. There exists an orthogonal matrix $P$ such that $P=(\bold p_1, \bold p_2, \cdots, \bold p_n)$ and
$A=P D P^t$ where
$$ D= \arr{\lambda_1 & 0 & \cdots & 0 \\ 0 & \lambda_2 &  \cdots & 0 \\ \vdots & \vdots & \ddots & \vdots \\
0 & 0 & \cdots & \lambda_n}. $$
Then
$$A-\lambda_i I =P D P^t-\lambda_i I= P(D-\lambda_i I) P^t=\sum_{j \neq i} (\lambda_j-\lambda_i) \mathbf p_j \mathbf p_j^t.$$
which equals 
$$\arr{p_{1,1} & \cdots & \widehat{p_{1,i}} & \cdots & p_{1,n} \\
p_{2,1} & \cdots & \widehat{p_{2,i}} & \cdots & p_{2, n} \\
\vdots & \ddots &  \vdots & \ddots & \vdots \\
p_{n, 1} & \cdots & \widehat{p_{n,i}} & \cdots & p_{n, n} } 
\arr{\lambda_1-\lambda_i  & \cdots & 0 & \cdots & 0 \\ 
\vdots & \ddots & \vdots & \ddots & \vdots \\
0 & \cdots & \widehat{\lambda_i-\lambda_i} & \cdots & 0 \\
\vdots & \ddots & \vdots & \ddots & \vdots \\
0 &  \cdots & 0 &  \cdots & \lambda_n-\lambda_i}
\arr{p_{1,1} & p_{2,1} & \cdots & p_{n,1} \\
\vdots & \vdots & \ddots & \vdots \\
\widehat{p_{1,i}} & \widehat{p_{2,i}} & \cdots & \widehat{p_{n,i}} \\
\vdots & \vdots &  \ddots & \vdots \\
p_{1,n} & p_{2,n} & \cdots & p_{n,n} }  .$$
Deleting the $m$-th row and $m$-th column, we obtain
$$ \arr{p_{1,1} & \cdots & \widehat{p_{1,i}} & \cdots & p_{1,n} \\
\vdots & \ddots &  \vdots & \ddots & \vdots \\
\widehat{p_{m,1}} & \cdots & \widehat{p_{m,i}} & \cdots & \widehat{p_{m, n}} \\
\vdots & \ddots &  \vdots & \ddots & \vdots \\
p_{n, 1} & \cdots & \widehat{p_{n,i}} & \cdots & p_{n, n} } 
\arr{\lambda_1-\lambda_i  & \cdots & 0 & \cdots & 0 \\ 
\vdots & \ddots & \vdots & \ddots & \vdots \\
0 & \cdots & \widehat{\lambda_i-\lambda_i} & \cdots & 0 \\
\vdots & \ddots & \vdots & \ddots & \vdots \\
0 &  \cdots & 0 &  \cdots & \lambda_n-\lambda_i}
\arr{p_{1,1} & \cdots & \widehat{p_{m,1}} & \cdots & p_{n,1} \\
\vdots & \ddots & \vdots & \ddots & \vdots \\
\widehat{p_{1,i}} & \cdots & \widehat{p_{m,i}} & \cdots & \widehat{p_{n,i}} \\
\vdots & \ddots &  \vdots & \ddots & \vdots \\
p_{1,n} & \cdots & \widehat{ p_{m,n}} & \cdots & p_{n,n} }  .$$
This is $A_m-\lambda_i I_{n-1}$. Notice that $P$ is orthogonal. Taking the determinant, we have
$$\det(A_m- \lambda_i I_{n-1})=p_{m,i}^2 \prod_{j \neq m} (\lambda_j-\lambda_i).$$
It follows that
$$p_{m,i}^2= \frac{\prod_{j=1}^{n-1} (\lambda_j(A_m)-\lambda_i)}{\prod_{j \neq m} (\lambda_j-\lambda_i)}.$$
Q.E.D.

\begin{cor}
Let $A$ and $B$ be two $n \times n$ real symmetric matrices. Suppose that
$eigen(A)= eigen(B)$ and $eigen(A_i)=eigen(B_i)$. Let $\lambda_i$ be a simple eigenvalue of $A$ and $B$. Let $\mathbf p_i=(p_{1,i}, p_{2,i}, \ldots, p_{n,i})^t$ be a unit vector in $eigen_{\lambda_i}(A)$ and 
$\bold q_i=(q_{1,i}, q_{2,i}, \ldots, q_{n,i})^t$ be a unit vector in $eigen_{\lambda_i}(B)$. Then
$$p_{j, i}^2=q_{j,i}^2 \ \forall j \in [1, n].$$
\end{cor}
\begin{cor}[Godsil-McKay, see Theorem 3.2, ~\cite{gm}]
Let $A$ and $B$ be two $n \times n$ real symmetric matrices. Suppose that
$A$ and $B$ are hypomorphic. Let $\lambda_i$ be a simple eigenvalue of $A$ and $B$. Let $\mathbf p_i=(p_{1,i}, p_{2,i}, \ldots, p_{n,i})^t$ be a unit vector in $eigen_{\lambda_i}(A)$ and 
$\bold q_i=(q_{1,i}, q_{2,i}, \ldots, q_{n,i})^t$ be a unit vector in $eigen_{\lambda_i}(B)$. Then
$$p_{j, i}^2=q_{j,i}^2 \ \forall j \in [1, n].$$
\end{cor}
\section{Eigenvalues and Eigenvectors under the perturbation of a rank one symmetric matrix}
Let $A$ be a $n \times n$ real symmetric matrix. Let $x$ be a $n$-dimensional row column vector. Let $M=x x^t$.
Now consider $A+ t  M$. We have
$$A+ t M= P D P^t + t M=P(D+ t P^{t} M P)P^t=P(D+ t P^{t} x x^t P)P^t.$$
Let $P^t x= q$. So $q_i=(\mathbf p_i, x)$ for each $i \in [1,n]$. Then
$$A+t \mathbf J= P( D+t q q^t) P^t.$$
Put $D(t)=D+ t q q^t$. 
\begin{lem}~\label{det} $\det(D+t q q^t-\lambda I) =\det(A-\lambda I) (1+\sum_{i} \frac{ t q_i^2}{\lambda_i -\lambda}).$
\end{lem}
Proof:  $\det(D-\lambda I+ t q q^t)$ can be written as a sum of products of $\lambda_i-\lambda$ and $q_i q_j$. For each $S $ a subset of $[1,n]$, combine the terms containing only $\prod_{i \in S} (\lambda_i -\lambda)$. Since the rank of $q q^t$ is one, only for $|S|=n, n-1$, the coefficients may be nonzero. We obtain
 $$det(D+t q q^t-\lambda I) =\prod_{i=1}^n (\lambda_i- \lambda)+ \sum_{i=1}^n t q_i^2 \prod_{j \neq i} (\lambda_i -\lambda).$$
 The Lemma follows. $\Box$\\
\\
Put $P_t(\lambda)=1+ \sum_{i} \frac{ t q_i^2}{\lambda_i -\lambda}$.
\begin{lem}~\label{fn} Fix $t<0$. Suppose that for each $i$, $\lambda_i$ is  a simple eigenvalue and $q_i \neq 0$. Then $P_t(\lambda)$ has exactly $n$ roots $(\mu_1, \mu_2, \cdots, \mu_n)$ satisfying a interlacing relation:
$$\lambda_1 > \mu_1 > \lambda_2 > \mu_2 > \cdots > \mu_{n-1} > \lambda_n > \mu_n.$$
\end{lem}
Proof: Clearly, $\frac{d P_t(\lambda)}{ d t}=\sum_{i} \frac{ t q_i^2}{(\lambda_i -\lambda)^2} < 0$. So $P_t(\lambda)$ is always decreasing. On the interval $(-\infty, \lambda_n)$, 
$P_t(-\infty)=1$ and $P_t(\lambda_n^-)=-\infty$. So $P_t(\lambda)$ has a unique root $\mu_n \in (-\infty, \lambda_n)$. Similar statement holds for each $(\lambda_{i-1}, \lambda_i)$. Q.E.D.
\begin{thm}~\label{main0}
 Fix $t < 0$. Let $l$ be the number of distinct eigenvalues satisfying
$(x, eigen_{\lambda}(A)) \neq 0$. Without loss of generalities, suppose that $A=P D P^t$ such that there exists a
$$S= \{ {{i_1}} > {{i_2}} > \cdots > {{i_l}} \}$$
 satisfying $ (x, \mathbf p_{i_j}) \neq 0 $ and  $(x, \mathbf p_i)=0$ for every $i \notin S $.  Then there exists 
$(\mu_1, \ldots, \mu_l)$ such that
$$\lambda_{i_1} > \mu_1 > \lambda_{i_2} > \mu_2 > \cdots > \lambda_{i_l} > \mu_l$$
and
$$eigen(A+ t M )=\{ \lambda_i(A) \mid i \notin S \} \cup \{\mu_1, \mu_2 \ldots, \mu_l \}.$$
Furthermore, $eigen_{\mu_j}(A+ t M )$ contains 
$$ \sum_{i \in S} \mathbf p_i \frac{q_i}{\lambda_{i}-\mu_j}.$$
\end{thm}
Here the index set $\{i_1, i_2, \cdots, i_l \}$ may not be unique.  I shall also point out a similar statement holds for $t>0$ with 
$$\mu_1 > \lambda_{i_1} > \mu_2 > \lambda_{i_2}  > \cdots > \mu_l >  \lambda_{i_l}.$$
Proof:  Since $ (x, eigen_{\lambda_{i_j}}(A)) \neq 0 $, $ q_{i_j} \neq 0$. For $i \notin S$, $q_i =0$. Notice 
$$P_t(\lambda)=1+\sum_{j=1}^{l} \frac{ t  q_{i_j}^2}{\lambda_{i_j}-\lambda}.$$
Applying  Lemma ~\ref{fn} to $S$, we obtain the roots of $P_t(\lambda)$
$\{ \mu_1, \mu_2, \ldots, \mu_l \}$ satisfying
$$\lambda_{i_1} > \mu_1 > \lambda_{i_2} > \mu_2 > \cdots > \lambda_{i_l} > \mu_l.$$
It follows that the roots of 
$\det(A+t M- \lambda I)=\det(D(t)-\lambda I)=P_t(\lambda) \prod_{i=1}^n (\lambda_i-\lambda)$
can be obtained from $eigen(A)$ be changing $\{ {\lambda_{i_1}} > {\lambda_{i_2}} > \cdots > {\lambda_{i_l}} \}$ to $\{\mu_1, \mu_2 \ldots, \mu_l \}$. Therefore,
$$ eigen(A+ t M)= \{ \lambda_i(A) \mid  i \notin S \} \cup \{\mu_1, \mu_2 \ldots, \mu_l \}.$$
For the sake of convenience, suppose that $\mu_i \notin eigen(A)$. Then 
$$ \sum_{i \in S}  \frac{q_i}{\lambda_{i}-\mu_j} \mathbf p_i=\sum_{i=1}^n  \frac{q_i}{\lambda_{i}-\mu_j} \mathbf p_i.$$
Here for $\lambda_i \notin S$, $q_i=0$.
Notice that
 $$ (A+ t M ) \sum_{i=1}^n  \frac{q_i}{\lambda_{i}-\mu_j} \mathbf p_i
= P(D+ t q q^t) P^t \sum_{i=1}^n  \frac{q_i}{\lambda_{i}-\mu_j} \mathbf p_i =P(D+t q q^t) \arr{\frac{q_1}{\lambda_1-\mu_j} \\ \vdots \\ \frac{q_n}{\lambda_n-\mu_j}},$$  
which equals
$$P \left( \arr{\frac{\lambda_1 q_1}{\lambda_1-\mu_j} \\ \vdots \\ \frac{\lambda_n q_n}{\lambda_n-\mu_j}}+
t \arr{q_1 \\ \vdots \\ q_n} \sum_{i=1}^n \frac{q_i^2}{\lambda_i-\mu_j} \right)= P \left( \arr{\frac{\lambda_1 q_1}{\lambda_1-\mu_j} \\ \vdots \\ \frac{\lambda_n q_n}{\lambda_n-\mu_j}}-
 \arr{q_1 \\ \vdots \\ q_n} \right)=P \arr{\frac{\mu_j q_1}{\lambda_1-\mu_j} \\ \vdots \\ \frac{\mu_j q_n}{\lambda_n-\mu_j}}.$$
We have obtained that $$(A+ t M) \sum_{i=1}^n  \frac{q_i}{\lambda_{i}-\mu_j} \mathbf p_i= \mu_j P \arr{\frac{ q_1}{\lambda_1-\mu_j} \\ \vdots \\ \frac{ q_n}{\lambda_n-\mu_j}}= \mu_j \sum_{i=1}^n  \frac{q_i}{\lambda_{i}-\mu_j} \mathbf p_i.$$
If $\mu_j \in eigen(A)$, we still have $(A+ t M ) \sum_{\lambda_i \in S} \frac{q_i}{\lambda_{i}-\mu_j} \mathbf p_i =
\sum_{i \in S} \frac{q_i}{\lambda_{i}-\mu_j} \mathbf p_i $. Therefore,
$$ \sum_{i \in S}  \frac{q_i}{\lambda_{i}-\mu_j} \mathbf p_i \in eigen_{\mu_j}(A+t M ).$$
Q.E.D.

\section{Reconstruction of Simple Eigenvectors not perpenticular to $\mathbf 1$}
Now let $M=\mathbf J=\mathbf 1 \mathbf 1^t$. Theorem ~\ref{main0} applies to $A+ t \bold J$ and $B+ t \bold J$.
\begin{thm}[Godsil-McKay, ~\cite{gm}] Let $B$ and $A$ be two real $ n \times n$ symmetric matrices.
Let $\Sigma$ be a hypomorphism such that $B=\Sigma(A)$.  Then there exists a subset $S \subseteq [1,n]$ such that $A=PDP^t$ and $B=UDU^t$ as in Theorem ~\ref{main0}. For $i \in S$, we have $\mathbf p_i=\mathbf u_i$ or $\mathbf p_i=-\mathbf u_i$. In particular, if $\lambda_i$ is a simple eigenvalue of $A$ and $(eigen_{\lambda_i}(A), \bold 1) \neq 0$, then $eigen_{\lambda_i}(A)=eigen_{\lambda_i}(B)$.
\end{thm}
Proof: $\bullet$ By Tutte's theorem, $eigen(A)=eigen(B)$. Let $A=PDP^t$ and $B=U D U^t$. Since $\det(A+ t \bold J -\lambda I)=\det(B+ t \bold J- \lambda I)$. By Lemma ~\ref{det}, 
$$\det(A-\lambda I) (1+\sum_{i} \frac{ t (\bold 1, \mathbf p_i)^2}{\lambda_i -\lambda})=\det(B-\lambda I)(1+\sum_{i} \frac{ t (\bold 1, \mathbf u_i)^2}{\lambda_i -\lambda}).$$
It follows that for every $\lambda_i$, $\sum_{\lambda_j=\lambda_i} (\bold 1, \mathbf p_j)^2=\sum_{\lambda_j=\lambda_i} (\bold 1, \mathbf u_j)^2$. Consequently, the $l$ for $A$ is the same as the $l$ for $B$. Let $S$ be as in Theorem ~\ref{main0} for both $A$ and $B$. Without loss of generality, suppose that $A=PDP^t$ and $B=UDU^t$ as in Theorem ~\ref{main0}. In particular,  for every $i \in [1,n]$, we have 
\begin{equation}~\label{1}
(\mathbf p_i, \bold 1)^2= (\mathbf u_i, \bold 1)^2.
\end{equation}
 $\bullet$ Let $T$ be as in Theorem ~\ref{main} for $A$ and $B$. Without loss of generality, suppose $T=(t_1, t_2) \subseteq \mathbb R^{-}$. Let $\mu_{l}(t)$ be the $\mu_l$ in Theorem ~\ref{main0} for $A$ and $B$. Notice that the lowest eigenvectors of $A+t \bold J$ and $B+ t \bold J$ are in ${\mathbb R^+}^n$ and they are not perpenticular to $\mathbf 1$. By Theorem ~\ref{main0}, $\mu_l(t)= \lambda_n(A+t \bold J)=\lambda_n(B+ t \bold J)$. By Theorem ~\ref{main},
$$eigen_{\mu_1(t)} (A+t \bold J)= eigen_{\mu_l(t)}(B+ t \bold J) \cong \mathbb R.$$
So $$\sum_{i \in S} \mathbf p_i \frac{(\mathbf p_i, \bold 1)}{\lambda_{i}-\mu_l(t)} // \sum_{i \in S} \mathbf u_i \frac{(\mathbf u_i, \bold 1)}{\lambda_{i}-\mu_l(t)}.$$
Since $\{ \mathbf p_i \}$ and $\{\mathbf u_i \}$ are orthogonal, by Equation ~\ref{1},
$$ \| \sum_{i \in S} \mathbf p_i \frac{(\mathbf p_i, \bold 1)}{\lambda_{i}-\mu_l(t)} \|^2 = \| \sum_{i \in S} \mathbf u_i \frac{(\mathbf u_i, \bold 1)}{\lambda_{i}-\mu_l(t)} \|^2.$$
It follows that for every $t \in T$,
$$\sum_{i \in S} \mathbf p_i \frac{(\mathbf p_i, \bold 1)}{\lambda_{i}-\mu_l(t)} = \pm \sum_{i \in S} \mathbf u_i \frac{(\mathbf u_i, \bold 1)}{\lambda_{i}-\mu_l(t)}.$$
$\bullet$ Recall that $-\frac{1}{t}= \sum_{i} \frac{ q_i^2}{\lambda_i -\mu_l(t) }$. Notice that the function $\rho \rightarrow \sum_{i} \frac{ q_i^2}{\lambda_i -\rho }$ is a continuous and one-to-one mapping from $(-\infty, \lambda_n)$ onto $(0, \infty)$. There exists a nonempty interval $T_0 \subseteq (-\infty, \lambda_n)$ such that if $\rho \in T_0$, then $ \sum_{i} \frac{ q_i^2}{\lambda_i -\rho} \in (-\frac{1}{t_1}, -\frac{1}{t_2})$. So every $\rho \in T_0$ is a $\mu_{l}(t)$ for some $t \in (t_1, t_2)$. It follow that  for every $\rho \in T_0$, $$\sum_{i \in S} \mathbf p_i \frac{(\mathbf p_i, \bold 1)}{\lambda_{i}-\rho} = \pm \sum_{i \in S} \mathbf u_i \frac{(\mathbf u_i, \bold 1)}{\lambda_{i}-\rho}.$$
Notice that both vectors are nonzero and depend continuously on $\rho$. Either,
$$\sum_{i \in S} \mathbf p_i \frac{(\mathbf p_i, \bold 1)}{\lambda_{i}-\rho} =  \sum_{i \in S} \mathbf u_i \frac{(\mathbf u_i, \bold 1)}{\lambda_{i}-\rho} \qquad \forall \ (\rho \in T_0);$$
or,
$$\sum_{i \in S} \mathbf p_i \frac{(\mathbf p_i, \bold 1)}{\lambda_{i}-\rho} =  - \sum_{i \in S} \mathbf u_i \frac{(\mathbf u_i, \bold 1)}{\lambda_{i}-\rho} \qquad \forall \ (\rho \in T_0);$$
$\bullet$. Notice that the functions $ \{ \rho \rightarrow  \frac{1}{\lambda_{i_{j}}-\rho} \}|_{i_j \in S}$ are linearly independent. For every $i \in S$, we have
$$\mathbf p_{i} (\mathbf p_i, \bold 1)= \pm \mathbf u_{i} (\mathbf u_i, \bold 1).$$
Because $\mathbf p_{i}$ and $\mathbf u_{i}$ are both unit vectors, $\mathbf p_{i}= \pm \mathbf u_{i}$. In particular, for every simple $\lambda_{i}$ with $(\mathbf p_{i}, \bold 1) \neq 0$ we have $eigen_{\lambda_i}(A)=eigen_{\lambda_i}(B)$. Q.E.D.
\begin{cor}
Let $B$ and $A$ be two real $ n \times n$ symmetric matrices.
Let $\Sigma$ be a hypomorphism such that $B=\Sigma(A)$. Let $\lambda_i$ be an eigenvalue of $A$ such that $(eigen_{\lambda_i}(A), \bold 1) \neq 0$. Then the orthogonal projection of $\bold 1$ onto $eigen_{\lambda_i}(A)$ equals the orthogonal projection of $\bold 1$ onto $eigen_{\lambda_i}(B)$.
\end{cor}
Proof: Notice that the projections are $\bold p_i(\bold p_i, \bold 1)$ and $\bold u_i(\bold u_i, \bold 1)$. Whether $\bold p_i= \bold u_i$ or $\bold p_i = - \bold u_i$, 
$$\bold p_i(\bold p_i, \bold 1)=  \bold u_i(\bold u_i, \bold 1).$$
Q.E.D.
\begin{conj} Let $A$ and $B$ be two hypomorphic matrices. Let $\lambda_i$ be a simple eigenvalue of $A$. Then there exists a permutation matrix $\tau$ such that
$\tau eigen_{\lambda_i}(A)= eigen_{\lambda_i}(B)$.
\end{conj}
This conjecture is apparently true if $eigen_{\lambda_i}(A)$ is not perpenticular to $\bold 1$.

\end{document}